# Modified Viterbi Algorithm Based Distribution System Restoration Strategy for Grid Resiliency

Chen Yuan, *Student Member, IEEE*, Mahesh S. Illindala, *Senior Member, IEEE*, and Amrit S. Khalsa, *Member, IEEE*

*Abstract*—This paper presents a novel modified Viterbi algorithm to identify the optimal distribution system restoration plan for improving the grid resiliency. In the proposed algorithm, the switching operations performed for system restoration are the states with the minimum bus voltage being seen as the cost metric for each state and the extent of load recovery as the observed event. When full load recovery is spotted, the dynamic programming algorithm stops, thereby giving the least number of switching pairs necessary for system restoration. Moreover, an improved flexible switching pair operation is employed to maintain the radial nature of distribution system. Several case studies are presented for verifying the performance of the proposed strategy. Multi-fault conditions are considered in testing the system restoration scheme on 33-bus and 69-bus distribution systems. Furthermore, the effects of integrating distributed energy resources and microgrid systems are analyzed.

## I. Introduction

POWER system blackouts are rare but extreme events in developed nations. Natural disasters and cascading failures due to overloads are the two main causes of large blackouts in North America. Weather events such as 2005 hurricane Katrina and 2012 hurricane Sandy are considered external sources [1], [2]. An overloading of power system could also trigger large outage following cascading failures as witnessed in the 2003 northeast blackout [3]. Hence, there is a critical need to incorporate resiliency into the power grid [4]–[6]. In radial distribution systems, when a blackout happens, the downstream customers also get disrupted during fault isolation. For estimating the cost of blackouts, the Lawrence Berkeley National Laboratory has published many reports [7]–[9] and a free web-based tool, the Interruption Cost Estimation (ICE) Calculator [10]. The principal aim for having a resilient system is fast recovery after an extreme event. Therefore, an efficient restoration strategy is necessary for quickly bringing back power to the un-faulted portions of grid network.

Faster detection of outages and plans for rapid restoration were envisioned in [11], [12] to build more resilient and secure power grids. Distribution system restoration strategies to recover electricity service to disrupted loads were presented in [13], [14]. The optimization of a restoration plan is a nonlinear problem with power flow operation limits and topology constraints. Several works were published to solve this problem, using heuristic methods [15]–[17], expert system [18], [19], multiagent system [20]–[22], fuzzy logic [23]–[25], neural network [26], mathematical programming [27]–[29], and graph theory [30], [31]. However, many focused on analysis of real power alone and neglected the reactive power. This is unacceptable because a reactive power deficit could lead to bus voltage drop, and even system collapse. Moreover, emphasis was primarily on fault at single location and finding its restoration plan. Very few works tested the restoration plan for worse cases involving multiple faults in different branches. The authors in [32] applied a spanning tree search algorithm to find the restoration plan and proposed a topology simplification scheme for large distribution system. They also considered the microgrids as virtual feeders in the distribution system. However, this paper does not give enough details on the microgrid architecture and capabilities. In [33]–[35], remote-controlled switches were used to maximize the restoration capability with minimum switch upgrade cost. But only single-fault conditions were considered. A dynamic programming based restoration plan was presented in [36]. It reduces the number of states to simplify the problem by grouping neighboring states and selecting the best one for each group. However, such a simplification affects the power flow analysis and compromises the restoration speed, thus violating the primary goal of quickly restoring the maximizing load.

In North America, approximately 3% of all distribution circuits have automated circuit reconfiguration (ACR). The existing ACR devices typically take 1-3 minutes to bring back power supply. By contrast, the remaining circuits have manual restoration that may take on an average 2-3 hours for being restored after the power outage is reported. The priority for restoration is given to critical infrastructures such as hospitals, 911 call centers, water treatment plants, police and fire stations, and then to largest pockets of customers. It should be noted that, according to the IEEE Guide for Electric Power Distribution Reliability Indices [37], outages less than 5 minutes are not counted as sustained interruptions. Therefore, such momentary outages are not considered in the reliability metrics like SAIDI, SAIFI, etc.

Manuscript received April 21, 2016; revised July 24, 2016; accepted September 7, 2016. Date of publication September 27, 2016; date of current version January 20, 2017. This work was supported in part by the Office of Naval Research under Award N00014-16-1-2753. Paper no.-TPWRD-00537-2016.
C. Yuan and M. S. Illindala are with the Department of Electrical and Computer Engineering, The Ohio State University, Columbus, OH 43210 USA (e-mail: yuan.261@osu.edu; millindala@ieee.org).
A. S. Khalsa is with the Dolan Technology Center, American Electric Power, Groveport, OH 43125 USA (e-mail: askhalsa@aep.com).
Color versions of one or more of the figures in this paper are available online at http://ieeexplore.ieee.org.
Digital Object Identifier 10.1109/TPWRD.2016.2613935

This paper proposes a novel modified Viterbi algorithm to identify the optimal restoration plan for power distribution systems with both real and reactive power flow analyses. The switching pair operation characterizes the state, whose cost metric is the minimum bus voltage within the distribution system. Besides, the load recovery rate is the observation information. In contrast to the conventional Viterbi algorithm, the observed event can be only estimated instead of being observed in advance. This estimated observed event also serves as a stop signal for the proposed modified Viterbi algorithm. Once the load is fully restored, the program is stopped, thus achieving the shortest switching operation sequence for system restoration. By this way, it ensures maximum load restoration in minimum number of switching operations. Furthermore, an improved flexible switching pair operation is proposed to maintain the distribution system's radial architecture. Various case studies are presented to verify the effectiveness of the proposed modified Viterbi algorithm for both single-fault and multi-fault conditions. The effectiveness in the presence of distributed energy resources (DERs) and microgrids is also demonstrated.

The rest of the paper is organized as follows. Section II presents a bi-level optimization problem formulation for distribution system restoration and the improved flexible switching pair operation. Section III illustrates the modified Viterbi algorithm proposed in this paper and the algorithm's innovative approach to minimize the switching operations in achieving maximal load restoration. Several case studies on benchmark 33-bus and 69-bus distribution systems with DERs and microgrids are presented in Section IV. Finally, the conclusion is given in Section V.

## II. PROBLEM FORMULATION

### A. Bi-Level Optimization Problem

The restoration plan has two hierarchical objectives—(i) to maximize load restoration, and then (ii) to realize it within the shortest time. The primary aim is to maximize the load restoration. To carry out the restoration plan, switching operations have to be undertaken either from remote control room or by dispatching field crews. If the number of switching operations could be minimized for achieving maximal load restoration, the restoration plan is optimized. Therefore, the strategy for post-fault distribution system restoration is formulated as an optimization problem with two hierarchical objectives, whose primary objective must be always satisfied and the secondary objective has the next preference.

Primary Objective

$$\text{Maximize} \sum_{i \in IB} S_i \quad (1)$$

Secondary Objective

$$\text{Minimize } N_{SW} \quad (2)$$

*Subject to*

$$S_i^p = P_i^p + jQ_i^p, \ i \in B, \ p \in X \quad (3)$$

$$P_i^p = |V_i^p| \sum_{j \in B} \sum_m |V_j^m| \cdot \left[ G_{ij}^{\text{pm}} \cos \theta_{ij}^{\text{pm}} + B_{ij}^{\text{pm}} \sin \theta_{ij}^{\text{pm}} \right] \quad (4)$$

$$Q_i^p = |V_i^p| \sum_{j \in B} \sum_m |V_j^m| \cdot \left[ G_{ij}^{\text{pm}} \sin \theta_{ij}^{\text{pm}} - B_{ij}^{\text{pm}} \cos \theta_{ij}^{\text{pm}} \right] \quad (5)$$

$$\left| I_{ij}^p \right| \leq \left| I_{ij}^{\max} \right|, \ i,j \in B, \ p \in X \quad (6)$$

$$\left| V_i^{\min} \right| \leq |V_i^p| \leq |V_i^{\max}|, \ i \in B, \ p \in X \quad (7)$$

$$\theta_i^{\min} \leq \theta_i^p \leq \theta_i^{\max}, \ i \in B, \ p \in X \quad (8)$$

$$\text{Maintaining radial network structure} \quad (9)$$

where
$B$    A set of buses in the distribution system
$IB$    A set of isolated buses in the distribution system
$X$    A set of phases a, b, c
$N_{SW}$    Number of switching pair operations
$S_i$    Total apparent power (MVA) of the load at bus $i$
$P_i$    Total real power (MW) of the load at bus $i$
$Q_i$    Total reactive power (MVAr) of the load at bus $i$
$V_i$    Phase voltage at bus $i$
$I_{ij}$    Phase current between buses $i$ and $j$
$\theta_i$    Phase voltage angle at bus $i$
$G_{ij}$    Real component of the $3 \times 3$ admittance matrix of branch between bus $i$ and bus $j$
$B_{ij}$    Reactive component of the $3 \times 3$ admittance matrix of branch between bus $i$ and bus $j$
$i, j$    Variable of bus number
$p, m$    Phase variable

As indicated in (1), the primary objective of the restoration plan is to guarantee the maximum load recovery. The secondary objective, given by (2), minimizes the switching pair operations during the system restoration process. However, the secondary objective of switching minimization can be only done after satisfying the primary objective. The power balance constraints in (3)–(5) show that the power injections of phase p at bus i should be equal to the load demand of phase p at bus i. The current and voltage limits are shown in (6) and (7). Inequality (8) expresses the voltage angle limit. An additional constraint for maintaining the distribution system's radial structure is shown in (9)–to make sure smooth operation of existing protection schemes in power distribution system. This requires a switching pair operation, which consists of closing a tie-line that is accompanied by opening of another. More details on avoiding mesh networks are presented in the next subsection.

Since maximum load recovery is the pursuit of the restoration plan, the primary objective in (1) can become a constraint of the secondary objective in (2). This is because the constraint must be satisfied first before optimizing the objective function. By this way, the maximum load restoration will be ensured no matter how the restoration is executed. Then the modified problem's objective is to minimize the number of switching pair operations with the maximization of the load restoration embedded within its constraints. The distribution system restoration strategy can be reformulated as a bi-level optimization problem:

Objective

$$\text{Minimize } N_{SW}$$

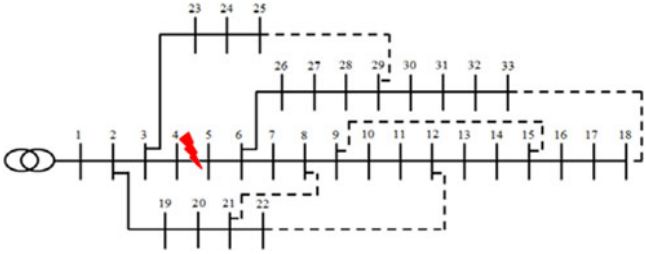

Fig. 1. 33-bus distribution system.

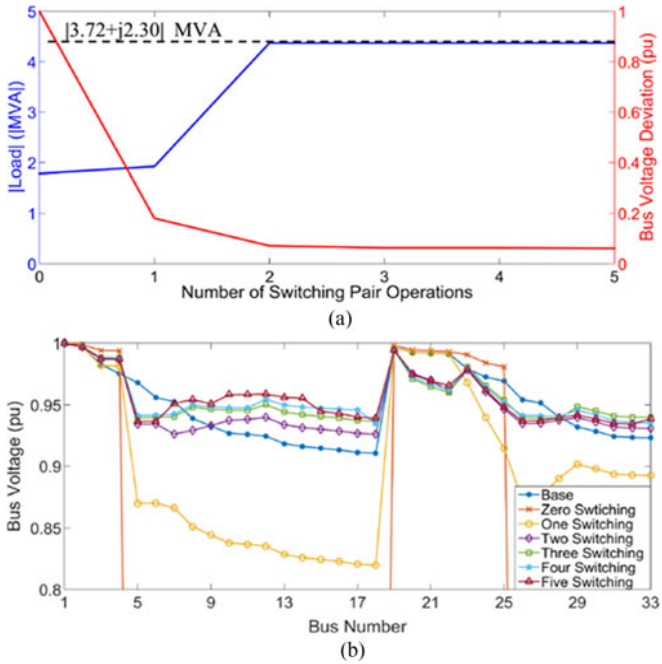

Fig. 2. Number of switching pair operations for (a) load recovery; (b) bus voltage distributions.

TABLE I
BEST SWITCHING PAIR OPERATIONS FOR THE CASE SHOWN IN FIG. 1

| Number of Switching Pair Operation | Best Switching Pair | Load Restoration | Maximum Bus Voltage Deviation |
|---|---|---|---|
| 0 | — | 0% | 100% |
| 1 | Close 25-29 | 5.3% | 18% |
| 2 | Close 25-29 & 12-22 Open 6-7 | 100% | 7.1% |
| 3 | Close 21-8 & 12-22 & 25-29 Open 11-12 & 28-29 | 100% | 6.3% |
| 4 | Close 21-8 & 12-22 & 25-29 & 18-33 Open 11-12 & 28-29 & 32-33 | 100% | 6.3% |
| 5 | Close 21-8 & 12-22 & 25-29 & 18-33 & 9-15 Open 9-10 & 14-15 & 17-18 & 28-29 | 100% | 6.3% |

In a severe case where an extreme event or a nature disaster can cause multiple faults, a more complex restoration plan is needed.

### B. Improved Flexible Switching Pair Operation

For a radial distribution system, after the fault isolation on the main feeder, the downstream loads are out-of-service. The power delivery should be restored quickly by closing a tie-line switch. However, this action may affect the distribution system's radial topology, thus needing switching pair operations. A switching pair is a group comprising one tie-line switch and one line segment. Since a mesh/loop is formed by closing a tie-line switch, to maintain the radial nature, this loop should be broken by opening a line segment within this loop. For each switching pair operation, the tie-line switch is easy to choose, but selecting a line segment to open is difficult. In previous works [29], [33], strategies with fixed switching pair and without fixed switching pair were considered. When the switching pair is fixed, operating one pair alone may not be sufficient for certain fault conditions. Severe challenges are observed in large systems. Therefore, the fixed switching pair is not a viable option.

It was explained in [29], [33] that if the switching pair is not fixed but flexible, the restoration plan can cover all fault conditions while maintaining radial network structure. However, when the search space is the whole system, it can take a long processing time to find the global optimal solution. Moreover, the complexity increases with growth of system size. Hence, in this paper, the flexible switching pair strategy limits the search space to within the loop—that was formed by closing the tie-line switch after isolating a fault—instead of the whole system. Once a tie-line switch is selected, finding the best switching pair is the same as finding the suitable line segment to open within the loop.

Hence, an improved flexible switching pair operation can be performed under two conditions, viz., (1) if the fault isolation is within the loop, no line segment needs to be opened when the tie-line closes, (2) otherwise, a line segment within the loop needs to open after the tie-line switch closes. For condition (1), fault isolation is indeed opening the line segment. Hence, the condition (1) can be treated as a special case of general switching

*Subject to*

$$\text{Maximize} \sum_{i \in IB} S_i$$

*Subject to*

(3)–(9)

The 33-bus distribution system, shown in Fig. 1, is used for illustration of this problem. Assume that a fault occurred on the distribution feeder between bus 4 and bus 5. Once the fault is isolated, the post-fault restoration begins. Fig. 2 shows the load recovery and bus voltage distributions. The load can be fully restored with just two switching pair operations. In addition, the bus voltage deviation is decreased to the minimum after three switching pair operations. All possible scenarios of switching operation are analyzed as explained further. For each set of switching pair operations, the best switching pair is chosen as the representative. Its selection is based on the load restoration level and the minimum bus voltage within the distribution system. The switching pairs used to obtain the results are presented in Table I. This is a simple case study for single-fault condition.

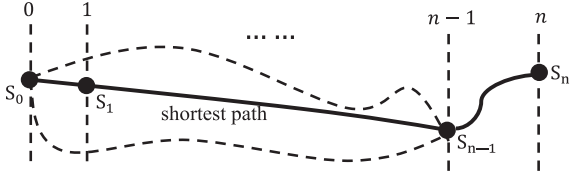

Fig. 3. Optimality principle of Viterbi algorithm.

pair operation covered in condition (2). Fig. 1 can be used for further explanation. After the fault isolation between bus 4 and bus 5, if the tie-line between buses 8 and 21 closes, there is no need to open any line segment. This is because the loop including buses 2-8 and 19-21 is broken by the fault isolation. However, if any bus voltage is lower than the minimum bus voltage limit, $|V_i^{\min}|$, another tie-line closure is needed. Under such circumstances, a distribution line segment may need to be opened to maintain a radial network structure. For example, if the tie-line between buses 25 and 29 is also closed, it forms a meshed loop with buses 2, 3, 6-8, 19-21 and 23-29. Therefore, a line segment between these buses should be opened.

## III. MODIFIED VITERBI ALGORITHM FOR DISTRIBUTION SYSTEM RESTORATION

### A. Viterbi Algorithm

The overall objective for distribution system restoration is to achieve maximum load restoration through minimum number of switching operations. To this end, the optimal restoration problem is a deterministic finite-state problem, equivalent to the shortest path problem. In practice, the shortest path can be easily constructed in sequence by forward dynamic programming via fast calculation of the shortest distance from the beginning through each stage. Viterbi algorithm is a forward dynamic programming approach for finding the most likely sequence of hidden states, called the Viterbi path, based on the observed events [38]. This algorithm lays out the states for each stage, uses a cost metric to evaluate all states, and determines the best sequence from the beginning stage to the final stage. Fig. 3 illustrates the optimality principle of Viterbi algorithm. Given the observation set $\bar{O} = \{O_0, O_1, \ldots, O_n\}$, the optimal sequence of states is $\bar{S} = \{S_0, S_1, \ldots, S_n\}$ with shortest path. If the previous $n-1$ segments of the shortest path to state $S_n$ are not the shortest path to state $S_{n-1}$ at the stage $n-1$, then there exists a shorter path to state $S_n$ [39]. So looking backward from each state at one stage, only the shortest path to this state is maintained.

### B. Modified Viterbi Algorithm for Distribution System Restoration

In contrast to the conventional Viterbi algorithm, a modified Viterbi algorithm is proposed for added benefits. The load recovery, as an observed event, is not known until the switching pair operation is executed. Each state signifies a switching pair operation, so the number of states of each stage depends on the number of switching pair operations, $C_m^k$, where $m$ is the number of tie-line switches and $k$ indicates the number of switching pair operation. Moreover, the number of stages is varying, and it depends on the restoration process. In other words, the ultimate number of switching pair operations depends on the extent of load recovery. If the load is fully recovered by satisfying all operational constraints, no more switching operations are needed. The corresponding state, instead of the sequence of states, is the preferred switching operation for restoration. Another distinction is in the cost metric. For the system restoration application, the minimum bus voltage is the cost index.

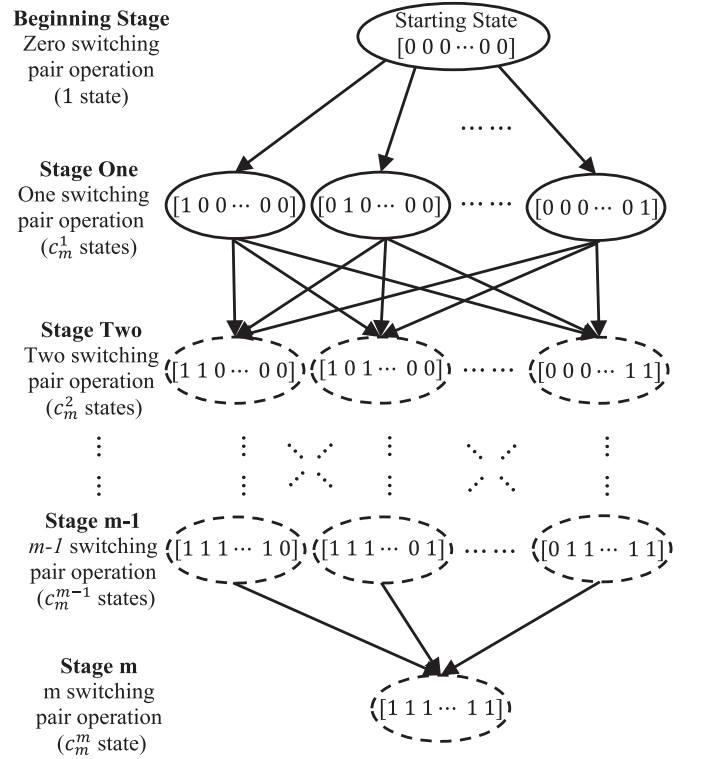

Fig. 4. State diagram of the modified Viterbi algorithm.

In the proposed modified Viterbi algorithm, a binary set is a state where each binary digit indicates the status of a tie-line switch status (0 stands for open and 1 means close). The number of binary digits denotes the total number of tie-line switches. This is because, in a radial distribution system, when a fault occurs, the downstream load is disconnected from the grid after the fault isolation. However, with the help of tie-line connection, power to the unserved load downstream can be restored. Since the tie-line switch is normally open, in the beginning stage, the starting state is a set of zeros. Thereafter, with the switching pair operations, the load recovery and bus voltage can be estimated by carrying out power flow analysis.

The state diagram of proposed modified Viterbi algorithm is illustrated in Fig. 4. To minimize the number of switching pair operations, the program ends once the load is fully recovered after meeting the constraints (3)–(9). However, the exact number of switching pair operations is not known until the modified Viterbi algorithm is applied. Therefore, the dashed lines in Fig. 4 show the potential states. The starting state is a set of zeros. At the first stage, one switching pair operation is implemented, so each state has only single one. At the second stage, two switching pair operations are implemented, with two ones in each state. In this sequence, at the kth stage, k switching pair operations are implemented. Hence, k ones are present in every state. If there exist m tie-line switches, the number of states in the first stage is designated as $C_m^1$. Likewise, it is $C_m^2$ in the second stage and $C_m^k$ in the kth stage. For each stage, the power flow analysis is conducted for every state—to check operational constraints (3)–(9). The program ends when the load is fully restored and all constraints are satisfied. The corresponding state is selected as the restoration scheme.

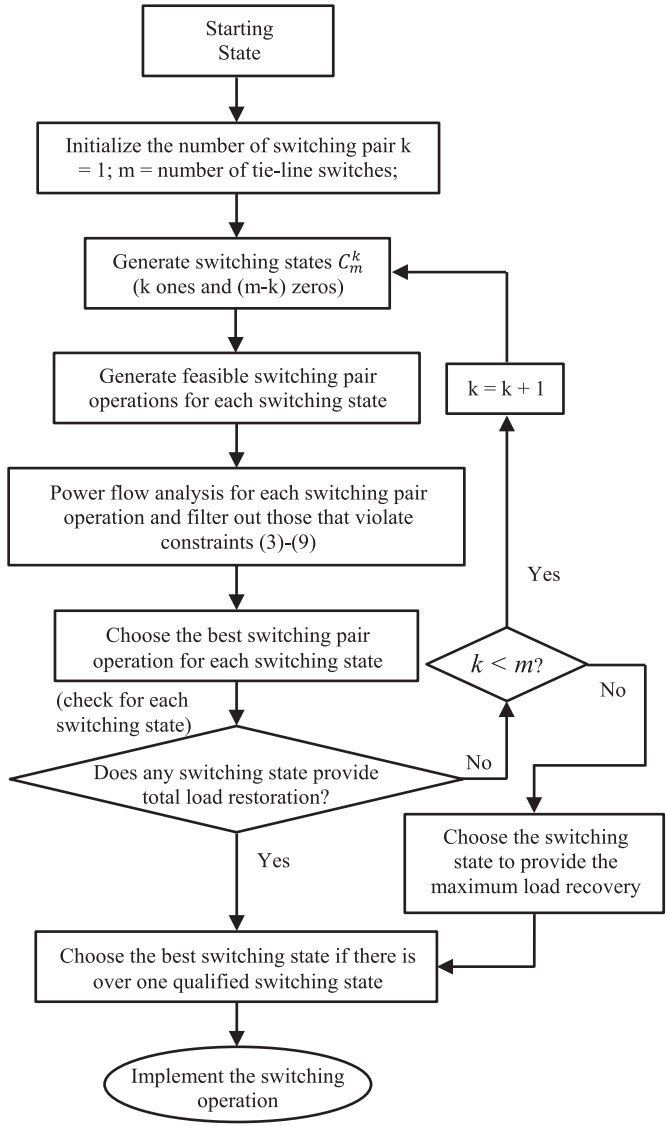

Fig. 5. Flowchart of the modified Viterbi algorithm based restoration scheme.

TABLE II
RESULTS OF RESTORATION FOR 33-BUS DISTRIBUTION SYSTEM

| Scenario | Fault Location | Switching Pair Operation | Load Restoration | Minimum In-Service Bus Voltage |
|---|---|---|---|---|
| 1 | 4-5 | Close 25-29 & 12-22 Open 6-7 | 100% | 0.9262 p.u. |
| 2 | 11-12 | Close 12-22 | 100% | 0.9322 p.u. |
| 3 | 4-5 & 27-28 | Close 25-29 & 8-21 | 100% | 0.9110 p.u. |
| 4 | 4-5 & 11-12 & 27-28 | Close 9-15 & 8-21 & 25-29 | 100% | 0.9191 p.u. |
| 5 | 5-6 & 6-7 & 6-26 | Close 25-29 & 12-22 | Load at bus 6 is unserved | 0.9262 p.u. |
| 6 | 13-14 & 15-16 & 6-26 & 26-27 | Close 9-15 & 18-33 & 25-29 | Load at bus 26 is unserved | 0.9190 p.u. |

The flowchart in Fig. 5 further describes the modified Viterbi algorithm based restoration scheme. At first, power flow analysis is carried out for all conditions of one switching pair operation. If the load can be fully restored, the switching pair operation is qualified. If over one switching pair operation could give the full load recovery, the one with highest *minimum bus voltage* is chosen for resiliency. In case multiple possibilities satisfy the above criterion, the one with lowest power loss is chosen. Otherwise, the choice is made randomly. If the load cannot be fully recovered by existing switching operations, it goes to next iteration by adding one more switching pair. At the end, if the load cannot be fully restored after closing all tie-line switches, the best switching pair operation is chosen to achieve the maximum load recovery. As shown in the flowchart, additional criteria are used, when necessary, to narrow the choices. When many switching options are available for maximal load recovery, the one with minimum switching is preferred. If no switching pair operation is implemented, it means no feasible switching pair operation satisfies (3)–(9). To restore maximal load, the load with the lowest bus voltage is cut off to relieve the system, and then the modified Viterbi algorithm is employed to find the best switching pair operation again.

## IV. CASE STUDIES

In this section, benchmark distribution systems are used to verify the proposed modified Viterbi algorithm. Both 33-bus and 69-bus distribution systems are utilized for testing on different system sizes. Furthermore, the 69-bus distribution system restoration is tested in the presence of distributed energy resources (DERs) and microgrid systems. The worst case is considered for each fault scenario to verify the effectiveness of the proposed algorithm. The distribution system load data, bus connections, line configuration and impedance are programmed into a MATLAB .m file. The modified Viterbi algorithm is also implemented in MATLAB software. Besides, the power flow analysis is conducted by using MATPOWER, which is a MATLAB power system simulation package.

### A. 33-Bus Distribution System

Fig. 1 shows a single-line diagram of the 33-bus distribution system [40], [41]. It is operating at 12.66 kV and supplying a total load of 3.72 MW and 2.30 MVAr. As shown in Fig. 1 by black dotted lines, five tie-lines are included to support post-fault system restoration. The allowable minimum voltage at any load bus is 0.9 p.u. and the maximum voltage is 1.05 p.u.

Six system restoration scenarios are presented in Table II. The first two are single-fault scenarios. In scenario 1, the system is restored with two switching pair operations, closing tie-line switches 25-29 and 12-22 and opening line-segment 6-7, and this results in minimum bus voltage is 0.9262 p.u. In scenario 2, the system can be fully restored only by closing tie-line switch 12-22. Scenario 3 has two faults, and so it needed closing tie-line switches 25-29 and 8-21. No extra line segment needed to be opened. Similarly, in scenario 4, with three fault isolations, only three tie-line switches are closed to restore the system. However, in scenarios 5 and 6, there exists isolated bus, when two neighboring faults happen. In these two cases, the load at the isolated bus cannot be restored. Therefore, the flowchart for maximizing restoration is implemented. Furthermore, the number of switching pairs is less than the number of faults—since two neighboring faults and isolated bus can be seen as one combined fault. In Fig. 6, the bus voltage distributions in different scenarios are presented and compared with each other. As seen in Fig. 6(b), the voltage of bus 6 is zero in scenario 5, indicating the bus isolation. The same condition applies to bus 26 in scenario 6.

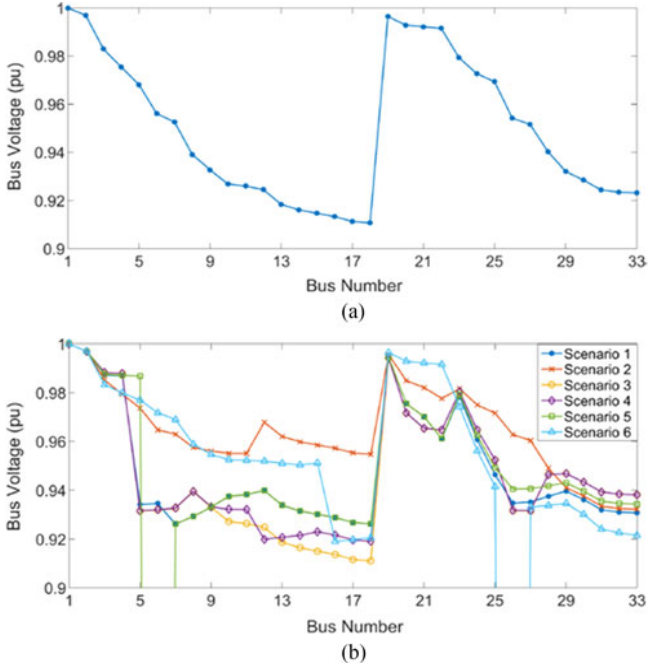

Fig. 6. Bus voltage distribution in (a) original 33-bus distribution system; (b) six scenarios.

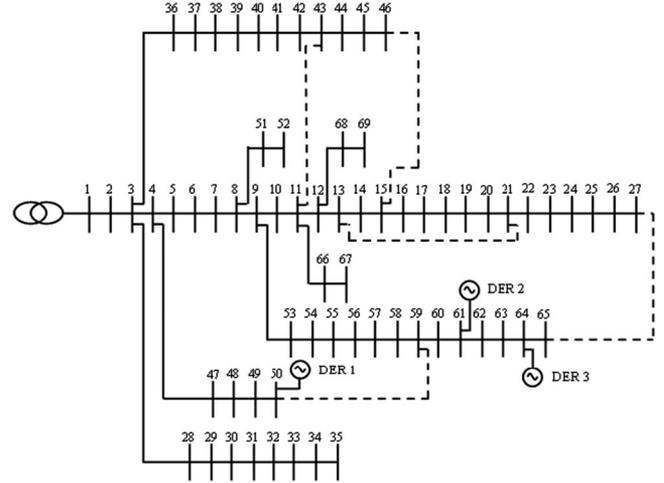

Fig. 7. 69-bus distribution system.

TABLE III
RESULTS OF RESTORATION FOR 69-BUS DISTRIBUTION SYSTEM

| Scenario | Fault Location | Switching Pair Operation | Load Restoration | Minimum In-Service Bus Voltage |
|---|---|---|---|---|
| 1 | 5-6 | Close 11-43 & 50-59 Open 55-56 (56-57 or 57-58 or 58-59) | 100% | 0.9349 p.u. |
| 2 | 18-19 | Close 13-21 | 100% | 0.9092 p.u. |
| 3 | 49-50 | Close 50-59 & 27-65 Open 61-62 | 100% | 0.9092 p.u. |
| 4 | 60-61 | Close 27-65 Open 61-62 | Load at bus 61 is unserved | 0.9407 p.u. |
| 5 | 42-43 & 43-44 | Close 11-43 & 15-46 | 100% | 0.9084 p.u. |
| 6 | 11-12 & 12-13 & 15-16 | Close 13-21 & 15-46 | Load at buses 12, 68 and 69 is unserved | 0.9142 p.u. |

### B. 69-Bus Distribution System Restoration With Penetration of DERs and Microgrid Systems

The 69-bus distribution system is shown in Fig. 7 [42], [43], including three DERs. It has 69 buses and 73 branches along with 5 normally open tie-lines. The nominal voltage is 12.66 kV, and the total load is 3.80 MW and 2.69 MVAr. The acceptable voltage range is 0.9 p.u. to 1.05 p.u. In this case study, the impact of DER integration in the 69-bus distribution system is also analyzed. DERs are installed at heavily loaded buses to reduce power losses and relieve the congestion on the distribution lines. In particular, they are assumed to be frequency-droop controlled voltage sources, and therefore the buses at which they are connected become PV buses. Placement of the DERs close to the critical loads will improve the system's resiliency for unexpected events in the distribution system.

Table III shows the six fault scenarios. The switching pair operation generated by the proposed Viterbi algorithm can achieve maximum load restoration. The corresponding restoration plans are presented in Table III. In scenario 1, the out-of-service load can be fully restored with two switching pair operations. In this condition, the four line segments between buses 55-59 can be equivalent, since there is no load on buses 55, 56, 57 and 58. There will be no power flow to these buses no matter which line segment is disconnected. In scenarios 2, 3 and 5, the loads are fully restored. However, in scenario 4, bus 61 is isolated because no feasible solution exists. Therefore the lowest voltage bus, 61, is disconnected, and then an optimal solution is found. In scenario 6, because of the neighboring faults, bus 12 is isolated, and the branch having buses 12, 68 and 69 is de-energized after fault isolation. Hence, the load at the three buses cannot be restored since no tie-line connects to this branch. Fig. 8 illustrates the bus voltage distributions for various scenarios. In scenarios 4

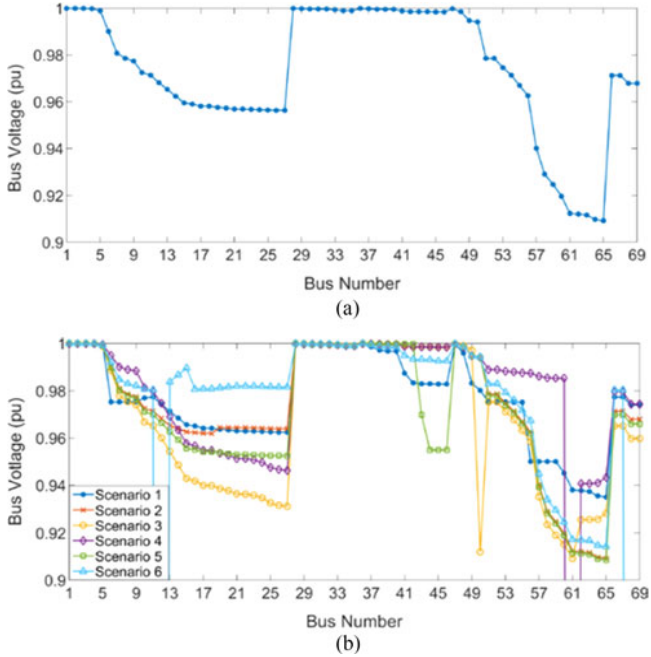

Fig. 8. Bus voltage distribution in (a) 69-bus distribution system; (b) six scenarios.

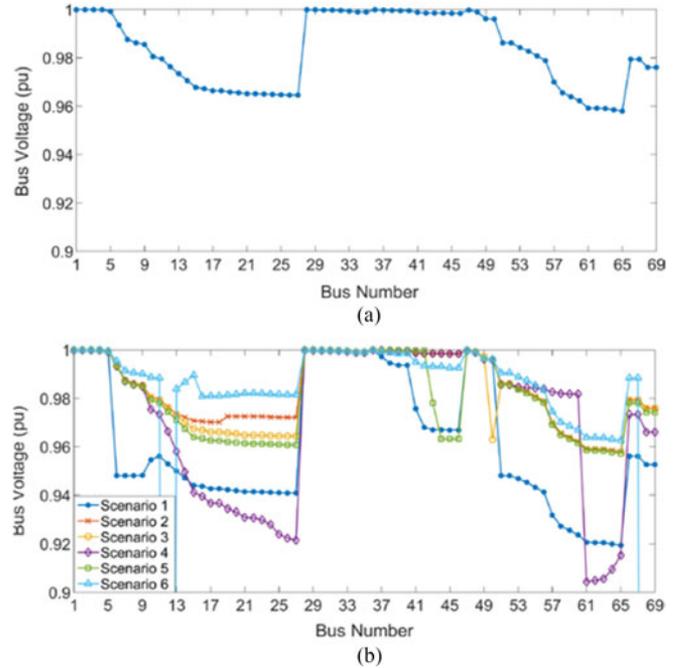

Fig. 9. Bus voltage distribution with DER integration in (a) 69-bus distribution system; (b) six scenarios.

TABLE IV
LOCATION AND GENERATION CAPACITY OF EACH DER

|  | DER 1 | DER 2 | DER 3 |
|---|---|---|---|
| Location | Bus 50 | Bus 61 | Bus 64 |
| Power Capacity | 0.5 MW | 1 MW | 0.3 MW |

TABLE V
RESULTS OF RESTORATION FOR 69-BUS DISTRIBUTION SYSTEM WITH DER INTEGRATION

| Scenario | Fault Location | Switching Pair Operation | Load Restoration | Minimum In-Service Bus Voltage |
|---|---|---|---|---|
| 1 | 5-6 | Close 11-43 | 100% | 0.9194 p.u. |
| 2 | 18-19 | Close 13-21 | 100% | 0.9580 p.u. |
| 3 | 49-50 | Close 50-59 | 100% | 0.9574 p.u. |
| 4 | 60-61 | Close 27-65 | 100% | 0.9043 p.u. |
| 5 | 42-43 & 43-44 | Close 11-43 & 15-46 | 100% | 0.9572 p.u. |
| 6 | 11-12 & 12-13 & 15-16 | Close 13-21 & 15-46 | Load at buses 12, 68 and 69 is unserved | 0.9626 p.u. |

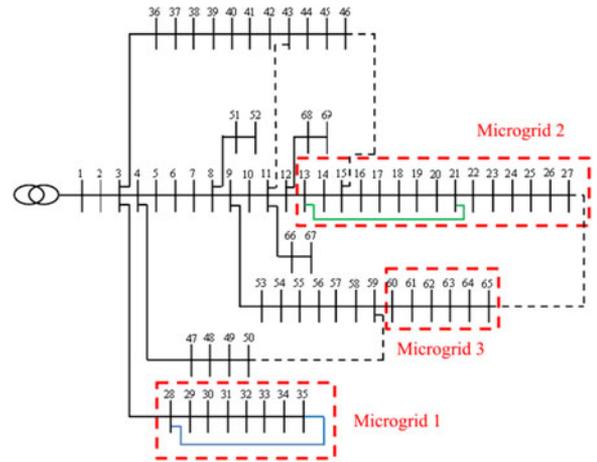

Fig. 10. 69-bus distribution system with three microgrid systems.

and 6, the isolated buses have zero voltage, indicated as less than 0.9 p.u. in Fig. 8(b).

Based on the load data in 69-bus system [44], three DERs are installed at buses 50, 61 and 64, as shown in Fig. 7. The three buses deliver 48.83% real power and 49.24% reactive power of the load. The location and generation capacity of each DER is displayed in Table IV.

To investigate the impact of DER integration on 69-bus distribution system restoration, analysis is carried out and results for the same fault scenarios as before are listed in Table V. However, in this case the switching pair operations are different. For scenarios 1–3, only one switching pair operation is needed. In scenario 4, bus 61 does not need to be isolated, because the DERs on bus 61 and bus 64 relieve the load burden. For scenarios 5 and 6, even though the switching pair operations are the same, the minimum in-service bus voltage is increased, thereby improving the system's reliability and resiliency after restoration. In Fig. 9, the bus voltage distributions in 69-bus system with DERs are presented. By comparing between Fig. 8 and Fig. 9, it is clear that the bus voltages are improved with DER integration.

Furthermore, the influence of microgrid systems is also studied. Fig. 10 has microgrid systems illustrated in red dashed regions within the 69-bus distribution system. In microgrid 1, a meshed topology is planned by adding an extra line, shown as the blue solid line between buses 28 and 35 (cf. Fig. 10). Microgrid 2 has a radial backbone, but forms a meshed loop by closing the tie-line between buses 13 and 21, shown with a green solid line in Fig. 10. In addition, a third microgrid (viz., Microgrid 3) with radial structure is considered. It is assumed that the three microgrids are self-sufficient and can deliver power to local loads without support from the main grid. For this reason, the three microgrid systems are equivalent to zero-load buses from the 69-bus distribution system side. Accordingly, the three microgrids can offset 52.92% real power and 52.65% reactive power of the total load demand in the distribution system. The generation capacity of each microgrid system is shown in Table VI. In post-fault restoration, three scenarios are possible: (1) If a fault occurs and gets isolated within the meshed loop of a microgrid, no restoration is needed because the meshed loop becomes radial structure after fault isolation and the load can be supplied without interruption. (2) If a fault occurs on the radial feeder of either microgrid 2 or 3, the power flow analysis should be carried out. In case any operational constraint is violated, a restoration plan is to be implemented. (3) If the fault happens outside the three microgrid systems, the post-fault restoration plan is required.

TABLE VI
GENERATION CAPACITY OF MICROGRIDS

|  |  | Microgrid 1 | Microgrid 2 | Microgrid 3 |
|---|---|---|---|---|
| Power Capacity | MW | 0.1 | 0.4 | 1.75 |
|  | MVAr | 0.07 | 0.24 | 1.12 |

TABLE VII
RESULTS OF RESTORATION FOR 69-BUS DISTRIBUTION SYSTEM WITH MICROGRID SYSTEMS

| Scenario | Fault Location | Switching Pair Operation | Load Restoration | Minimum In-Service Bus Voltage |
|---|---|---|---|---|
| 1 | 5-6 | Close 11-43 | 100% | 0.9762 p.u. |
| 2 | 18-19 | — | 100% | 0.9894 p.u. |
| 3 | 49-50 | Close 50-59 | 100% | 0.9696 p.u. |
| 4 | 60-61 | — | 100% | 0.9894 p.u. |
| 5 | 42-43 & 43-44 | Close 11-43 & 15-46 | 100% | 0.9860 p.u. |
| 6 | 11-12 & 12-13 & 15-16 | — | Load isolated and unserved at buses 12, 68 and 69 | 0.9920 p.u. |

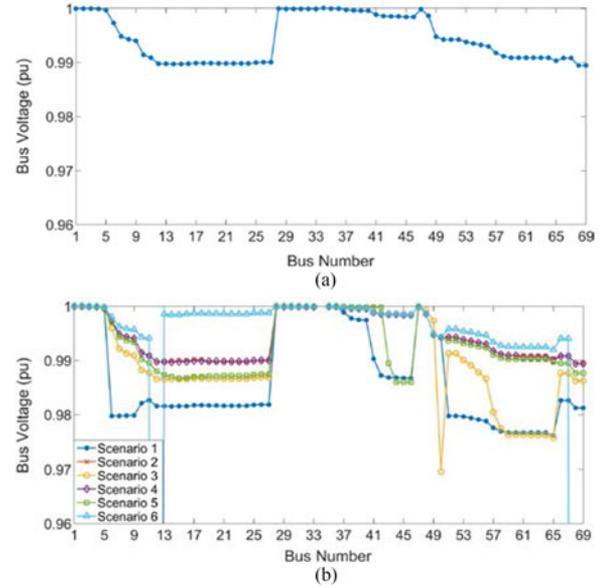

Fig. 11. Bus voltage distribution in the presence of microgrid systems for (a) 69-bus distribution system; (b) six scenarios.

Table VII shows the results of distribution system restoration with microgrid systems. For scenarios 1, 3 and 5, although the switching pair operations are the same as the system with DERs, the minimum bus voltages are much higher than before. For scenarios 2, 4, 6, no switching action is needed since the faults happen within the microgrid systems. Besides, no load is out of service in scenarios 2 and 4, but loads at buses 12, 68 and 69 are not restored in scenario 6. However, this is not possible since no tie-line is connecting to the branch through buses 12, 68 and 69. Fig. 11 illustrates the bus voltage distributions in 69-bus distribution system with the three microgrid systems. When compared with Fig. 8 and Fig. 9, Fig. 11 shows much improvement for the bus voltages. It is also noticed that the bus voltage distribution in scenario 2 is very close to that in scenario 4. This is because the faults happen within microgrid systems and the loads are covered by local DERs.

## V. CONCLUSION

For improving the resiliency of power grid, a novel modified Viterbi algorithm based distribution system restoration strategy was proposed. The strategy was to formulate it as a bi-level optimization problem to maximize load recovery with least number of switching pair operations following fault isolation. An improved flexible switching pair operation was used to maintain the system's radial architecture. Various case studies were presented to demonstrate the performance of the proposed restoration strategy on both 33-bus and 69-bus distribution systems. Both single fault and multi-fault conditions were considered. The power flow analysis and modified Viterbi algorithm have been implemented in MATLAB to obtain optimal restoration planning. Furthermore, the effectiveness of proposed restoration scheme in the presence of DERs and microgrids was also demonstrated.

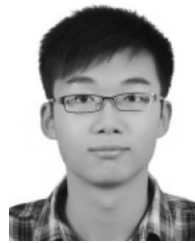
**Chen Yuan** (S'13) received the B.S. degree in electrical engineering from Wuhan University, Wuhan, China, in 2012 and is currently pursuing the Ph.D. degree in electrical and computer engineering from The Ohio State University, Columbus, OH, USA.

His current research interests include the control and protection of microgrids, distributed energy resources, and energy management systems.

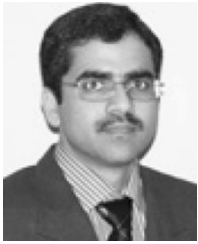

**Mahesh S. Illindala** (S'01–M'06–SM'11) received the B.Tech. degree in electrical engineering from the National Institute of Technology, Calicut, India, in 1995, the M.Sc.(Engg.) degree in electrical engineering from the Indian Institute of Science, Bangalore, India, in 1999, and the Ph.D. degree in electrical engineering from the University of Wisconsin, Madison, WI, USA, in 2005.

From 2005 to 2011, he was with Caterpillar R&D. Since 2011, he has been an Assistant Professor in the Department of Electrical and Computer Engineering, The Ohio State University, Columbus, OH, USA. His research interests include microgrids, distributed energy resources, electrical energy conversion and storage, power system applications of multiagent systems, as well as protective relaying and advanced electric drive transportation systems.

Prof. Illindala received the ONR Young Investigator Award in 2016.

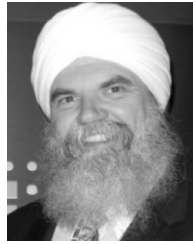

**Amrit Khalsa** (M'91) received the B.S. and M.S. degrees in electrical engineering from The Ohio State University, Columbus, OH, USA in 1984 and 2011, respectively.

From 1984 to 1987, he was a Teaching Assistant in the Department of Electrical Engineering, The Ohio State University. Since 1988, he has held various engineering and management positions with American Electric Power (AEP), Columbus, in the areas of system protection and control, underground systems engineering, distribution system planning, geospatial information systems, and advanced technology testing. He is currently a Staff Engineer at AEP's Dolan Technology Center. He is an author/coauthor of papers in areas including microgrids and electric transportation. His research interests include microgrids, distributed energy resources, and smart grid technology from the distribution system down to the customer.

Mr. Khalsa is a Registered Professional Engineer in the State of Ohio. He received the EPRI Technology Transfer Award in 2012 for contributions to developing ANSI Standard CEA-2045, a modular communication interface for demand response.